\documentclass{article}
\usepackage{amssymb,amsmath,amsthm, graphicx, epsfig}

\newtheorem{lemma}{Lemma}

\begin{document}

\title{The 2-adic valuation of plane partitions and totally symmetric plane partitions}
\author{William J. Keith}

\maketitle

\abstract{This paper confirms a conjecture of Amdeberhan and Moll that the power of 2 dividing the number of plane partitions in an n-cube is greater than the power of 2 dividing the number of totally symmetric plane partitions in the same cube when n is even, and less when n is odd.}

\section{Introduction}

In \cite{AmMoll}, Tewodros Amdeberhan and Victor H. Moll discussed the 2-adic valuation of the number of alternating sign matrices, plane partitions, and many of the distinguished symmetric subsets of plane partitions, such as totally symmetric self-complementary plane partitions.  They conjectured (Conjecture 3.1 in \cite{AmMoll}) that 

\[ v_2(PP_{2n}) > v_2(TSPP_{2n}) \, \text{ and } \, v_2(PP_{2n+1}) < v_2(TSPP_{2n+1}) \]

\noindent where $PP_n$ is the number of plane partitions in the $n \times n \times n$ cube, $TSPP_{n}$ is the number of totally symmetric plane partitions in the $n$-cube, i.e. plane partitions symmetric under any permutation of the axes, and $v_2(n)$ is the highest power of 2 dividing $n$.  The main theorem of this paper is that the conjecture holds.  In fact, considerably more seems to be true, the data for which inspired the proof strategy.

\subsection{Definitions}

A \emph{plane partition} $\pi$ of an integer $x$ is an array $\pi_{i,j}$ of nonnegative integers in $i,j \geq 1$ such that $\sum \pi_{i,j} = x$, which is nonincreasing in rows and columns, that is, $\pi_{i,j} \geq \pi_{i+1,j}$ and $\pi_{i,j} \geq \pi_{i,j+1}$ for all $i,j$.  This condition means that $\pi$ has finite support.  The points $\{(i,j,k) \vert 0 < k \leq \pi_{i,j}\}$ constitute the three-dimensional Young diagram of $\pi$.  We say $\pi$ is in the $n$-cube if $\pi_{i,j} \leq n$ for all $i,j$ and $\pi_{i,j} > 0 \Rightarrow i \leq n \, \text{ and } \, j \leq n$.  We denote the number of plane partitions in the $n$-cube by $PP_n$.  It has the formula

\[ PP_n = \prod_{i,j,k=1}^n \frac{i+j+k-1}{i+j+k-2} \, \text{.}\]

Plane partitions for which the Young diagram is invariant under various involutions of the cube have similarly compact formulas.  The set of plane partitions which are invariant under exchange of the $i$ and $j$ axes are called \emph{symmetric plane partitions}; their number in the $n$-cube is denoted $SPP_n$ and has formula \[ SPP_n = \prod_{j=1}^n \prod_{i=j}^n \frac{i+j+n-1}{i+j-1} \, \text{.} \]  The plane partitions invariant under \emph{any} permutations of the three axes are called \emph{totally symmetric}, and are counted by the formula \[ TSPP_n = \prod_{1\leq i \leq j \leq k \leq n} \frac{i+j+k-1}{i+j+k-2} \, \text{.} \]  Plane partitions invariant under permutations of the axes, as well as complementation -- the operation of taking all lattice points in the $n \times n \times n$ cube not in the Young diagram of the partition, and exchanging corners to make this a new diagram -- are called \emph{totally symmetric self-complementary} partitions, and they can only appear in a cube of even size 2n.  In such a cube, their counting formula is \[ TSSCPP_{2n} = \prod_{1\leq i \leq j \leq n} \frac{i+j+n-1}{i+j+i-1} \, \text{ .}\]

We denote the highest power of 2 dividing $n$ by $v_2(n)$, i.e. $v_2(n)=k$ if $2^k \mid n$ but $2^{k+1} \nmid n$.  We use $s_2(n)$ to denote the number of 1s in the binary expansion of $n$; this function has the properties that $s_2(2j+1) = 1 + s_2(j)$, $s_2(2j) = s_2(j)$ and $s_2(3m+1) \leq 2 s_2(m) + 1$.  These two quantities are related by a formula of Legendre, $v_2(m!) = m - s_2(m)$, which will repeatedly be useful in the following.

\section{Proof of the Main Result}

Two points to notice are that, first, the formula for $TSPP_n$ is a subset of the factors in the formula for $PP_n$, and second, a great deal of cancellation occurs in both formulas.  Triples $(i,j,k)$ from $PP_n$ are rejected for $TSPP_n$ if $i>j$, regardless of the relation of $k$ to $i$ and $j$, or if $i \leq j$ but $j > k$, whether $i>k$ or not.

Thus:

\begin{multline*} PP_n=TSPP_n \cdot \prod_{j=1}^n \prod_{i=j+1}^n \prod_{k = 1}^n \frac{i+j+k-1}{i+j+k-2} \cdot \prod_{i=1}^n \prod_{j = i}^n \prod_{k = 1}^{j-1} \frac{i+j+k-1}{i+j+k-2} \\
= TSPP_n \cdot \prod_{i=1}^n \prod_{j=i+1}^n \prod_{k = 1}^n \frac{i+j+k-1}{i+j+k-2} \cdot \prod_{i=1}^n \prod_{j = i}^n \prod_{k = 1}^{j-1} \frac{i+j+k-1}{i+j+k-2}
\end{multline*}

\noindent where for convenience we have switched the $i$ and $j$ entries in the middle product; the factors are symmetric in $i$ and $j$, so this introduces no problems.

Many of the factors in the numerator and denominator cancel.  In the middle term, $\prod_{j=1}^n \prod_{i=j+1}^n \prod_{k = 1}^n \frac{i+j+k-1}{i+j+k-2}$, most factors arising from a triple $(i,j,k)$ in the numerator cancel with the factor arising from $(i,j,k+1)$ in the denominator.  The exceptions are triples with $k=n$ in the numerator, and those with $k=1$ in the denominator.  Thus, this term becomes

\[ \prod_{i=1}^n \prod_{j=i+1}^n \frac{i+j+n-1}{i+j-1} =  SPP_n \cdot \prod_{i=1}^n \frac{2i-1}{2i+n-1} \, \text{.} \]

That is, this term is exactly the formula for $SPP_n$ except that it is missing the boundary where $i=j$.

Notice that $v_2(2i-1)=0$.  If $n$ is even, $2i+n-1$ is odd, and so $v_2 (2i+n-1)= 0$.  If $n$ is odd, $v_2(2i+n-1)>0$.  Over the whole product, we can add in intervening odd factors to create a factorial without changing the 2-adic valuation of the product, to get 

\begin{multline*} v_2( \prod_{i=1}^n 2i+n-1) = v_2( (3n-1)!) - v_2 ((n-1)!) \\
= 3n-1 - s_2(3n-1) - (n-1) + s_2(n-1) \\
= 2n - s_2(3n-1) + s_2(n-1) \, \text{.} \end{multline*}

In the third factor, $\prod_{i=1}^n \prod_{j=i}^n \prod_{k = 1}^{j-1} \frac{i+j+k-1}{i+j+k-2}$, a triple $(i,j,k)$ in the numerator cancels with $(i+1,j,k)$ in the denominator unless $i=j$ in the numerator (if $i=n$, $i=j$ automatically), or $i=1$ in the denominator.  We are left with

\begin{multline*} \prod_{i=1}^n \prod_{j=i}^n \prod_{k = 1}^{j-1} \frac{i+j+k-1}{i+j+k-2}  = \prod_{j=1}^n \prod_{k = 1}^{j-1} \frac{2j+k-1}{j+k-1}  \\
= \prod_{j=1}^n \frac{(2j)(2j+1)\cdots(3j-2)}{j(j+1)\cdots(2j-2)} \, \text{.}
 \end{multline*}
 
The 2-adic valuation of this product is
 
 \begin{multline*}
v_2 (\prod_{j=1}^n \frac{(2j)(2j+1)\cdots(3j-2)}{j(j+1)\cdots(2j-2)} ) \\
= \sum_{j=1}^n v_2((3j-2)!) - v_2((2j-1)!) - v_2((2j-2)!) + v_2((j-1)!) \\
= \sum_{j=1}^n 3j-2 - s_2(3j-2) - (2j-1) + s_2(2j-1) - (2j-2) + s_2(2j-2) + j-1 - s_2(j-1) \\
= \sum_{j=1}^n s_2(2j-1) + s_2(2j-2) - s_2(3j-2) - s_2(j-1) \\
= \sum_{j=1}^n s_2(j-1) + 1 + s_2(j-1) - s_2(3j-2) - s_2(j-1) \\
= \sum_{j=1}^n s_2(j-1) + 1 - s_2(3j-2) = n + \sum_{j=1}^n s_2(j-1) - s_2(3j-2) \, \text{.}
 \end{multline*}
 
 The last analysis we can perform before we must split the proof into the cases of $n$ even and odd is to determine $v_2(SPP_n)$:
 
\begin{multline*} v_2(SPP_n) = v_2 \left( \prod_{j=1}^n \frac{(2j+n-1)(2j+n)\cdots(j+2n-1)}{(2j-1)(2j)\cdots(j+n-1)} \right) \\
= v_2 \left( \frac{(3n-1)}{(2n-1)} \cdot \frac{(3n-2)(3n-3)}{(2n-2)(2n-3)} \cdot \frac{(3n-3)(3n-4)(3n-5)}{(2n-3)(2n-4)(2n-5)} \right. \times \\
\left.  \cdots \times \frac{(2n)(2n-1)\cdots(n+1)}{(n)(n-1)\cdots(1)} \right) \\
= v_2 \left( \frac{(3n-1)!}{(2n-1)!} \cdot \frac{(3n-3)!}{(2n-2)!} \cdot \frac{(3n-5)!}{(2n-3)!} \cdots \frac{(n+1)!}{n!} \right) \\ / \left( \frac{(2n-1)!}{(n-1)!} \cdot \frac{(2n-3)!}{(n-2)!} \cdots \frac{1!}{0!} ) \right) \\
= \sum_{j=1}^n 3n+1-2j - s_2(3n+1-2j) + n-j -s_2(n-j) \\
- (2n-j) + s_2(2n-j) - (2n+1-2j) + s_2(2n+1-2j) \\
= \sum_{j=1}^n s_2(2n-j) + s_2(2n+1-2j) - s_2(3n+1-2j) - s_2(n-j) \\
= n + \sum_{j=1}^n s_2(2n-j) - s_2(3n+1-2j) \, \text{.}
\end{multline*}
 
 Putting the pieces together so far, we have: if $n$ is even,
 
 \[ v_2(PP_n) = v_2(TSPP_n) + 2n + \sum_{j=1}^n s_2(2n-j) - s_2(3n+1-2j) + s_2(j-1) - s_2(3j-2) \]
 
\noindent and if $n$ is odd, 
 
 \begin{multline*} v_2(PP_n) = v_2(TSPP_n) + s_2(3n-1) - s_2(n-1) \\ + \sum_{j=1}^n s_2(2n-j) - s_2(3n+1-2j) + s_2(j-1) - s_2(3j-2) \, \text{.} \end{multline*}

We now examine the two cases separately.

\textbf{Case 1: Even.}  When $n$ is even, the claim is equivalent to showing that for even $n$, 

\begin{equation}\label{EvenCase} 2n + \sum_{j=1}^n s_2(2n-j) + s_2(j-1) -s_2(3j-2) - s_2(3n+1-2j) > 0\, \text{ .}
\end{equation}

Begin by noticing that the factors for $TSSCPP_{2n}$ are a subset of the factors counted in this expression.  We subtract the 2-adic valuation of $TSSCPP_{2n}$, which, since this is an integer, is nonnegative; the resulting expression will be simpler.

\begin{multline*}TSSCPP_{2n} = \prod_{1 \leq i \leq j \leq n} \frac{i+j+n-1}{i+j+i-1} \\
= \frac{(3n-1)}{(3n-1)} \cdot \frac{(3n-3)(3n-2)}{(3n-4)(3n-3)} \cdot \frac{(3n-5)(3n-4)(3n-3)}{(3n-7)(3n-6)(3n-5)} \cdots \frac{(n+1)(n+2)..(2n)}{(2)(3)..(n+1)} \\
= \left[ \frac{(3n-1)!}{(2n-1)!} \cdot \frac{(3n-3)!}{(2n-2)!} \cdot \frac{(3n-5)!}{(2n-3)!} \cdots \frac{(n+1)!}{(n)!} \right] \\ / \left[ \frac{(3n-1)!}{(3n-2)!} \cdot \frac{(3n-3)!}{(3n-5)!} \cdot \frac{(3n-5)!}{(3n-8)!} \cdots \frac{(n+1)!}{(1)!} \right]
\end{multline*}

Thus

\begin{multline*}v_2 (TSSCPP_{2n}) = v_2 \left( \frac{(3n-2)! (3n-5)! .. (1)!}{(2n-1)! (2n-2)! .. (n)!} \right) \\
= \sum_{j=1}^n 3j-2 - s_2(3j-2) - (j+n-1) + s_2(j+n-1) \\
= \sum_{j=1}^n 2j-1-n - s_2(3j-2) + s_2(j+n-1) \\
= \sum_{j=1}^n s_2(j+n-1) - s_2(3j-2) \, \text{ .}
\end{multline*}

Subtracting this from the left-hand side of equation (\ref{EvenCase}), we get

\begin{multline}
2n + \sum_{j=1}^n s_2(2n-j) - s_2(3n+1-2j) + s_2(j-1) - s_2(3j-2) - s_2(j+n-1) + s_2(3j-2) \\
= 2n + \sum_{j=1}^n s_2(j-1) - s_2(3n+1-2j) = n + \sum_{j=1}^n s_2(j-1) - s_2((3n/2)-j) \\
= n + \sum_{j=1}^n s_2(j-1) - s_2(j-1+(n/2)) = n + \sum_{j=0}^{n-1} s_2(j) - s_2(j+(n/2)) \\
= 2k + \sum_{j=0}^{2k-1} s_2(j) - s_2(j+k) = 2k + \sum_{j=0}^{k-1} s_2(j) - \sum_{j=0}^{k-1} s_2(j+2k)
\end{multline}

\noindent where we have set $k = \frac{n}{2}$.

Let $S_E(k) = 2k + \sum_{j=0}^{k-1} s_2(j) - \sum_{j=0}^{k-1} s_2(j+2k)$.  We wish to show that this expression is strictly positive -- that, if the total number of binary digits for all numbers from 0 to $k-1$ is compared to that for the interval of the same length moved up $2k$, the excess of binary digits in the latter over the former is no more than $2k$.

\begin{lemma}\label{DigitLemma} For $k>0$, we have $\sum_{j=0}^{k-1} s_2(j+2k) - \sum_{j=0}^{k-1} s_2(j) < 2k$.
\end{lemma}

\noindent \textbf{Proof: } Let $k = 2^{i_1} + \dots + 2^{i_r}$, $i_1 < \dots < i_r$.  We will concern ourselves with three contributions: 1s in the binary expansion of the larger values that appear in columns beyond $i_r$, the difference in the number of 1s appearing within the interval in column $i_r$, and the differences in the number of 1s appearing in all columns below $i_r$.

\begin{figure}[h]
\begin{center}
\begin{tabular}{c c}
\begin{tabular}{c|ccc | c | cc|}
 & 0 & 1 & 2 & $3$ & $4$ & 5 \\
 \hline 0 & 0 & 0 & 0 & 0 & 0 & 0 \\
 1 & 1 & 0 & 0 & 0 & 0 & 0  \\
 2 & 0 & 1 & 0 & 0 & 0 & 0  \\
 3 & 1 & 1 & 0 & 0 & 0 & 0  \\
 4 & 0 & 0 & 1 & 0 & 0 & 0  \\
 5 & 1 & 0 & 1 & 0 & 0 & 0  \\
 6 & 0 & 1 & 1 & 0 & 0 & 0  \\
 7 & 1 & 1 & 1 & 0 & 0 & 0  \\
\cline{2-5} 8 & 0 & 0 & 0 & 1 & 0 & 0  \\
 9 & 1 & 0 & 0 & 1 & 0 & 0  \\
 10 & 0 & 1 & 0 & 1 & 0 & 0  \\
 \hline 
\end{tabular}

&

\begin{tabular}{c|ccc | c | cc|}
 & 0 & 1 & 2 & $3$ & $4$ & 5 \\
 \hline 22 & 0 & 1 & 1 & 0 & 1 & 0 \\
 23 & 1 & 1 & 1 & 0 & 1 & 0  \\
 24 & 0 & 0 & 0 & 1 & 1 & 0  \\
 25 & 1 & 0 & 0 & 1 & 1 & 0  \\
 26 & 0 & 1 & 0 & 1 & 1 & 0  \\
 27 & 1 & 1 & 0 & 1 & 1 & 0  \\
 28 & 0 & 0 & 1 & 1 & 1 & 0  \\
 29 & 1 & 0 & 1 & 1 & 1 & 0  \\
\cline{2-5} 30 & 0 & 1 & 1 & 1 & 1 & 0  \\
 31 & 1 & 1 & 1 & 1 & 1 & 0  \\
 32 & 0 & 0 & 0 & 0 & 0 & 1  \\
 \hline 
\end{tabular}
\end{tabular}
\caption{The intervals of interest for $k=11$.}
\end{center}
\end{figure}

Since the smallest value in the interval $[2k,3k-1]$ is $2k$, there is exactly one 1 in column $i_r+1$ for 2k, but no 1s in any higher column.  It is never the case that columns $i_r+1$ and $i_r+2$ both contain a 1 for numbers in $[2k,3k-1]$, since $2^{i_r+1} + 2^{i_r+2} \geq 3k+3$.  Hence columns beyond $i_r$ contain exactly one 1 for each number in the interval $[2k,3k-1]$, contributing $k$ to the sum.  Thus we must show that less than $k$ additional 1s are added in the columns up to $i_r$.

The digits in column $i$ have period $2^{i+1}$, consisting of a string of $2^i$ 0s followed by $2^i$ 1s.  As long as $k \neq 2^{i_r}$, column $i_r$ does not complete a full period, so it consists of a string of 0s followed by a string of 1s of length $2^{i_1} + \dots + 2^{i_{r-1}}$.  (If $k = 2^{i_r}$, it is simple to notice that the frame of interest consists of even multiples of the periods of all columns of interest, and so the higher frame can consist only of rotations of these; so that the number of 1s counted within the higher frame cannot differ at all.  Hence we will assume $k \neq 2^{i_r}$ in the sequel.)

Consider $k^{-} = 2^{i_1} + \dots + 2^{i_{r-1}} \leq 2^{i_r} - 1$.  The digits in column $i_r$ will be rotated forward by $2k^{-}$, since adding $2 \cdot 2^{i_r}$ only moves the frame forward an even period in column $i_r$.  The difference $add(i_r)$ between the number of 1s in the resulting column, and those that existed in the original, will be a fraction of $2^{i_r}$ given by a piecewise linear function of the fraction of $2^{i_r}$ represented by $k^{-}$:

$$ \frac{add(i_r)}{2^{i_r}} = \left\{ \begin{matrix} 
2 \left( \frac{k^{-}}{2^{i_r}} \right) & \quad & \left( \frac{k^{-}}{2^{i_r}} \right) < \frac{1}{3} \\
1 - \left( \frac{k^{-}}{2^{i_r}} \right) & \quad & \frac{1}{3} < \left( \frac{k^{-}}{2^{i_r}} \right) \leq \frac{1}{2} \\
2 - 3 \left( \frac{k^{-}}{2^{i_r}} \right) & \quad & \frac{1}{2} < \left( \frac{k^{-}}{2^{i_r}} \right) < \frac{2}{3} \\
0 & \quad & \frac{2}{3} < \left( \frac{k^{-}}{2^{i_r}} \right) < 1 \end{matrix} \right.$$

\begin{center}\begin{tabular}{c}
\includegraphics[scale=0.75]{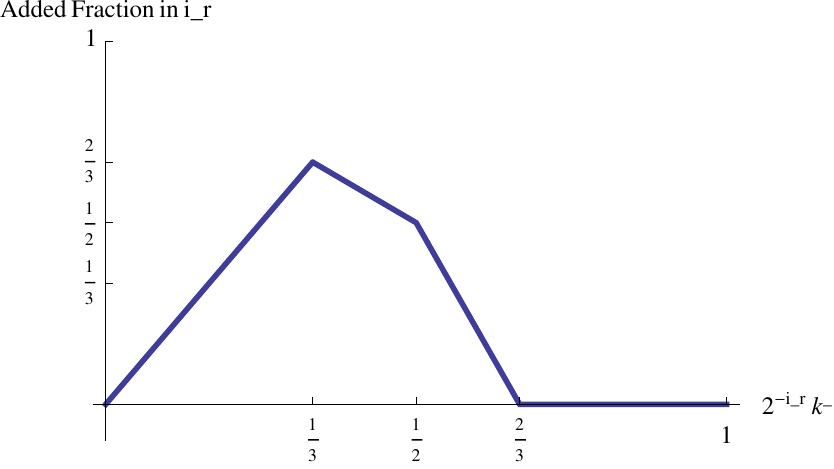}
\end{tabular}\end{center}

Therefore, at most, the additional 1s in column $i_r$ contribute $\frac{2}{3} 2^{i_r}$ to the difference in the digit sums, while the digit 1 in column $i_r$ contributes $2^{i_r}$ to the value of $k$.  If we are attempting to add $k$ 1s, this leaves a deficit of $\frac{1}{3} 2^{i_r}$ which, we will see, cannot be made up in lower columns.

In lower columns, the column within the frame of interest may or may not end with a 1.  Therefore, in looking at column b, whether $2^b$ is a digit in $k$ will also concern us.  In general, let $k_b = \sum_{i_j \leq b} 2^{i_j}$.  Then the additional number of 1s in column $b$ is given by the piecewise linear function

$$ \frac{add(b)}{2^{b}} = \left\{ \begin{matrix} 
0 & \quad & 0 \leq  \left( \frac{k_b}{2^{b+1}} \right) < \frac{1}{6} \\
2 \left( \frac{k_b}{2^{b+1}} \right) - \frac{1}{3} & \quad & \frac{1}{6} \leq  \left( \frac{k_b}{2^{b+1}} \right) < \frac{1}{3} \\
1 - 2 \left( \frac{k_b}{2^{b+1}} \right) & \quad & \frac{1}{3} \leq  \left( \frac{k_b}{2^{b+1}} \right) < \frac{1}{2} \\
-2 + 4 \left( \frac{k_b}{2^{b+1}} \right) & \quad & \frac{1}{2} \leq  \left( \frac{k_b}{2^{b+1}} \right) < \frac{2}{3} \\
\frac{10}{3} - 4 \left( \frac{k_b}{2^{b+1}} \right) & \quad & \frac{2}{3} \leq  \left( \frac{k_b}{2^{b+1}} \right) < \frac{5}{6} \\
0 & \quad & \frac{5}{6} \leq  \left( \frac{k_b}{2^{b+1}} \right) < 1
 \end{matrix} \right.$$

\begin{center}\begin{tabular}{c}
\includegraphics[scale=0.75]{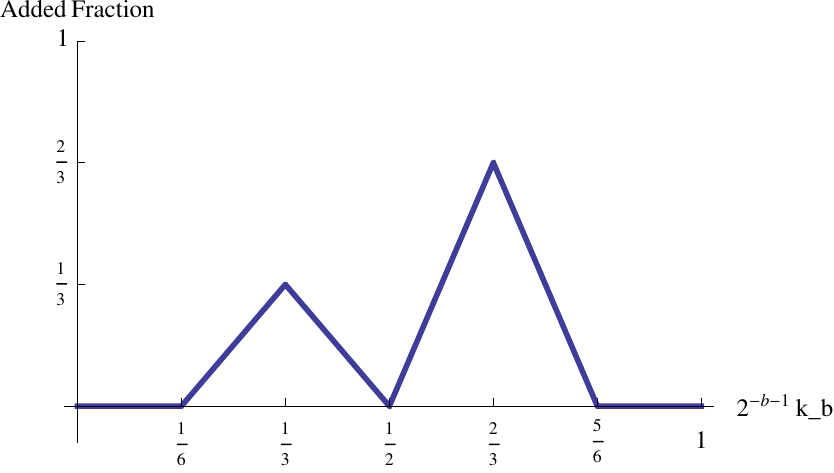}
\end{tabular}\end{center}

Anywhere $\frac{k_b}{2^{b+1}} \geq \frac{1}{2}$, the digit in $k$ in column $b$ must be a 1, contributing $2^b$ to the value of $k$.  In such a case, the extra 1s number less than the addition to the value of $k$, and so the deficit between $k$ and the number of added 1s increases.  If $\frac{k_b}{2^{b+1}} < \frac{1}{2}$, the most that can be added is $\frac{1}{3} 2^b$.  Supposing this to be the case for as many $b$ as possible, the sum still must be at most $\frac{1}{3} \left(1 + 2 + \cdots + 2^{i_r - 1} \right) < \frac{1}{3} 2^{i_r}$.  Thus each column either increases the deficit, or is insufficient a contribution to close the gap.

Hence fewer than $k$ 1s are added in the columns up to and including $i_r$, and exactly $k$ 1s are added in the columns beyond, so the additional 1s number less than $2k$.  The lemma, and thus this case of the theorem, is proved. $\Box$

\textbf{Case 2: Odd. }  For $n$ odd,

\begin{multline*}
\sum_{j=1}^n s_2(2n-j) - s_2(3n+1-2j) + s_2(j-1) - s_2(3j-2) - v_2(TSCPP_{2n}) \\
= \sum_{j=1}^n s_2(2n-j) - s_2(3n+1-2j) + s_2(j-1) - s_2(3j-2) - (\sum_{j=1}^n s_2(j+n-1) - s_2(3j-2)) \\
= \sum_{j=1}^n s_2(j-1) - s_2(3n+1-2j) = \sum_{j=0}^{n-1} s_2(j) - s_2(\frac{n+1}{2}+j) \\
= \sum_{j=0}^{\frac{n-1}{2}} s_2(j) - s_2(n+j) = \sum_{j=0}^{k-1} s_2(j) - s_2(j+2k+1) \, \text{ for } \, k = \frac{n-1}{2} \\
= s_2(n) - s_2(3(n-1)/2+1) + \sum_{j=0}^{k-1} s_2(j) - s_2(j+2k) \\
= 1 + s_2(n-1) - s_2(3n-1) + \sum_{j=0}^{k-1} s_2(j) - s_2(j+2k)
\end{multline*}

and so

\begin{multline*}
s_2(3n-1) - s_2(n-1) + \sum_{j=1}^n s_2(2n-j) - s_2(3n+1-2j) + s_2(j-1) - s_2(3j-2) \\
= v_2(TSCPP_{2n}) + 1 + \sum_{j=0}^{k-1} s_2(j) - s_2(j+2k) \, \text{ .}
\end{multline*}

By the arguments in the previous section, $\sum_{j=0}^{k-1} s_2(j) - s_2(j+2k)$ is negative, and at least $k$.  But $v_2(TSCPP_{2n})$ is known (\cite{MollSun}) to be much less than $\frac{n-1}{2}$: it reaches a maximum of the closest integer to $\frac{2^n}{3}$ on the interval $\left[\frac{2}{3} 2^n , \frac{4}{3} 2^n\right]$, and is less than half this value on this interval outside $\left[\frac{5}{6} 2^n , \frac{7}{6} 2^n\right]$.  Thus this sum is always negative, and so the conjecture is a theorem. $\Box$

\section{Relative Sizes}

These proof methods are limited in what they can say about the relative sizes of $v_2(PP_n)$ and $v_2(TSPP_n)$.  The data on the differences calculated, however, suggest two very interesting conjectures.

For even index, the difference $v_2(PP_{2n}) - v_2(TSPP_{2n})$ appears to be very nearly 5 times the value $v_2(TSSCPP_{2n})$.

\begin{center}\begin{tabular}{cc}
\includegraphics[scale=0.6]{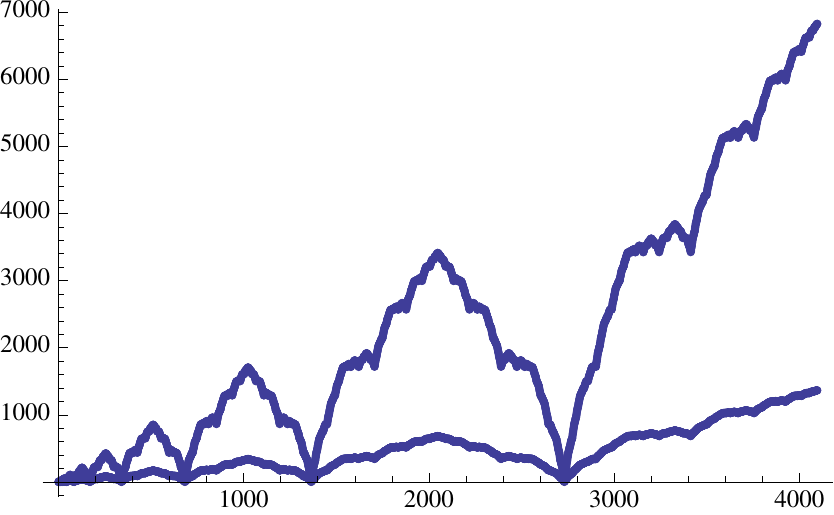} & \includegraphics[scale=0.6]{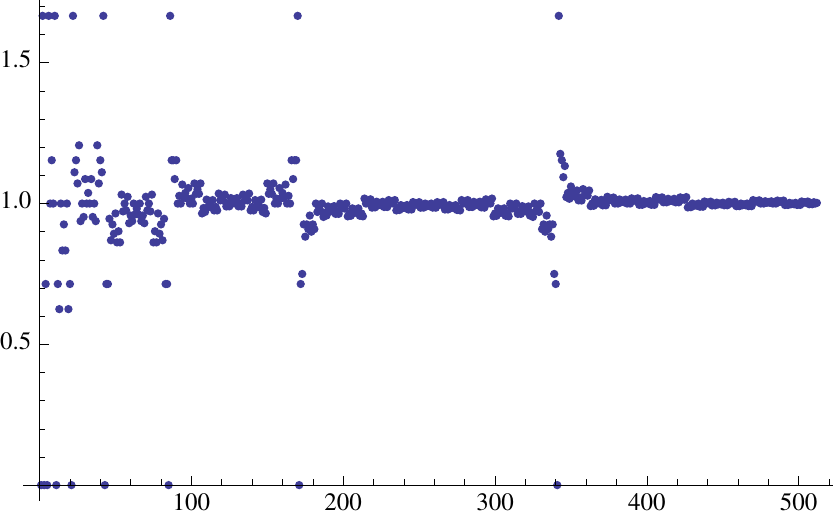}
\end{tabular}\end{center}

Here we have graphed, on the left, the values $v_2(PP_{2n}) - v_2(TSPP_{2n})$ and $v_2(TSSCPP_{2n})$.  On the right, we have the ratio $\frac{5 v_2(TSSCPP_{2n})}{v_2(PP_{2n}) - v_2(TSPP_{2n})}$.  The ratio is eccentric near the very low values of $v_2(TSSCPP_{2n})$, but for the most part of the interval is very close to 1.  The maximum ratios appear to be exactly 5/3.

For odd index, the difference $v_2(PP_{2n-1}) - v_2(TSPP_{2n-1})$ decreases steadily, its minima lying on the line $-3n$.  In the right-hand graph below we have diplayed $\frac{5 v_2(TSSCPP_{2n})}{v_2(PP_{2n-1}) - v_2(TSPP_{2n-1})+3n}$.  The maxima of this ratio are not constant, but again for the most part of the interval the ratio appears to be close to exact.

\begin{center}\begin{tabular}{cc}
\includegraphics[scale=0.6]{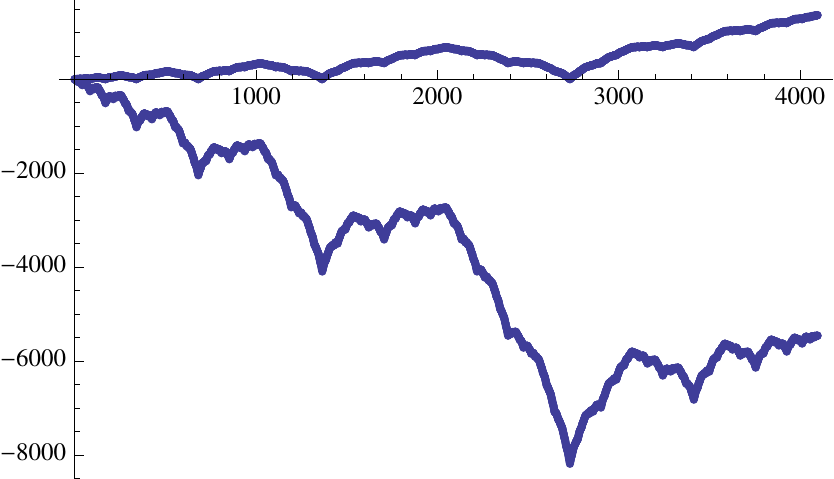} & \includegraphics[scale=0.6]{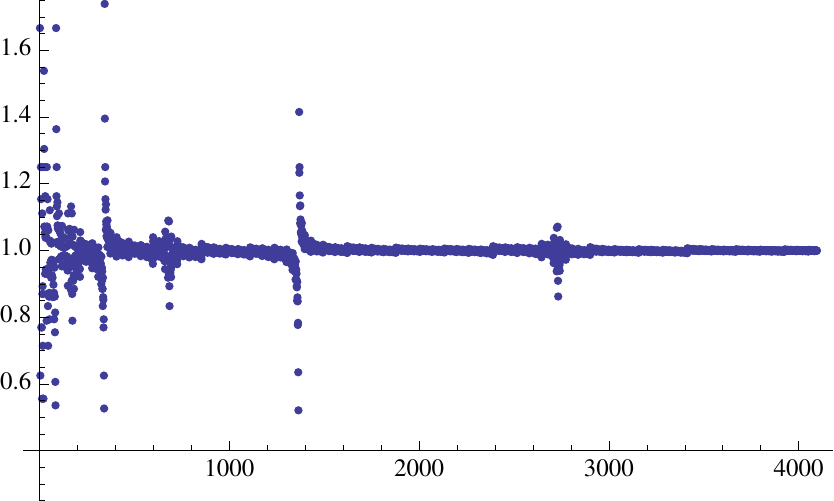}
\end{tabular}\end{center}

It is certainly conceivable that these conjectures could be made more precise, and it seems likely that the proof will employ the fact that the ordered triples comprising $TSPP_n$ are just 1/6 the cube of triples comprising $PP_n$ in their respective formula.  However, we offer the conjectures to the reader in their present form.

\end{document}